\documentclass[a4wide]{article}

 \usepackage{authblk}
\usepackage{float}
\usepackage{amsmath,amssymb,array} % estetics
 \usepackage{stmaryrd} 
\usepackage{graphicx,esvect}

\newcommand{\bfu}{{\boldsymbol u}}
\newcommand{\bfv}{{\boldsymbol v}}

\newcommand{\bff}{{\boldsymbol f}}
\newcommand{\bfn}{{\boldsymbol n}}

\newcommand{\bft}{{\boldsymbol t}}
\newcommand{\bfI}{{\boldsymbol I}}

\newcommand{\bfS}{{\boldsymbol S}}
\newcommand{\bfC}{{\boldsymbol C}}

\newcommand{\bftau}{{\boldsymbol\tau}}
\newcommand{\bfsig}{{\boldsymbol\sigma}}

\newcommand{\bfeps}{{\boldsymbol\varepsilon}}

\newcommand{\jump}[1]{\llbracket{#1}\rrbracket}

\newcommand{\IR}{\mathbb{R}}
\numberwithin{equation}{section}
\newtheorem{remark}{Remark}[section]
\date{}

\begin{document}
\title{Nitsche's Finite Element Method for Model Coupling in Elasticity}
% \date{\today}
\author[$\dagger$]{Peter Hansbo}
\author[$\ddagger$]{Mats G. Larson}
\affil[$\dagger$]{\footnotesize\it  Department of Mechanical Engineering, J\"onk\"oping University, SE-55111 J\"onk\"oping, Sweden}
\affil[$\ddagger$]{\footnotesize\it Department of Mathematics and Mathematical Statistics, Ume{\aa}~University, SE-901\,87 Ume{\aa}, Sweden}

\maketitle

\begin{abstract}
We develop a Nitsche finite element method for a model of Euler--Bernoulli beams
with axial stiffness embedded in a two--dimensional elastic bulk domain. The beams have their own displacement fields, and the elastic subdomains created by the beam network are triangulated independently and are coupled to the beams
weakly by use of Nitsche's method in the framework of hybridization. 
\end{abstract}
%\keywords{Nitsche's method, model coupling, hybrid method, interface stiffness}

\section{Introduction} 

In this paper we continue our work on coupling of elastic models\cite{BuHaLa17,CeHaLa16,Ha05,HaLaLa17}. 
Unlike previous coupling models \cite{Ha05,NgKeClBo13}, we take as our starting point
the hybridized approach of Burman et al.\cite{BuElHaLaLa19}, where an auxiliary interface variable is introduced in the Nitsche framework. The hybridized formulation conveniently supports solution and preconditioning based on substructuring 
where the bulk variables are eliminated resulting in a system for the hybrid variable. Furthermore, the hybrid variable 
may be used to model interfaces with mechanical properties such as bending and membrane stiffness. We consider in particular embedded interfaces made up by beams-trusses embedded in a two dimensional elastic membrane with both strong and cohesive coupling between the beam-truss interface and the membrane. Using the hybridized Nitsche framework we easily derive a weak formulation, which directly leads 
to a finite element method by replacing the function spaces with conforming finite dimensional finite element spaces. 
The cohesive formulation is designed in such a way that we may let the stiffness in the coupling tend to infinity without loss  of stability or convergence. We focus our attention on fitted meshes, that are not  required to match on 
the interface, but the approach can directly be extended to cut finite element formulations.

The outline of the paper is as follows: In Section 2 we introduce the hybridized formulation for an elastic interface problem, in Section 3 we consider interfaces with bending and membrane stiffness as well as strong and cohesive coupling to the elastic problem. In Section 4 we present some numerical examples illustrating the method and we observe  optimal order convergence properties.

\section{The Elastic Interface Problem}

We begin by extending the hybridized Nitsche method proposed by Burman et al. \cite{BuElHaLaLa19} to the case of linearized
elasticity. Let $\Omega$ denote a bounded domain in $\IR^2$. For ease of presentation, we consider only the case where
$\Omega$ is divided into two non-overlapping subdomains
$\Omega_{1}$ and $\Omega_{2}$, $\Omega_{1}\cup  \Omega_{2}$, with
interface $\Gamma=\overline\Omega_{1}\cap \overline \Omega_{2}$; the case of several domains is a straightforward extension.
We further assume that the subdomains are polygonal
so that $\Gamma$ is piecewise linear. In each subdomain, we assume plane stress linearized elasticity with homogeneous Dirichlet boundary conditions, i.e.,
we seek displacement fields $(\bfu_1,\bfu_2,\bfu_\Gamma)$ that are zero on $\partial\Omega$. The problem takes the form: find 
$\bfu_i:\Omega_i \rightarrow \IR^2$ and $\bfu_\Gamma :\Gamma\rightarrow \IR^2$
such that 
\begin{alignat}{3}
\bfsig(\bfu_i)   & = 2 \mu \bfeps(\bfu_i)+\lambda ~\nabla\cdot\bfu_i\bfI& \qquad &\text{in $\Omega_i$}
\\
-\bfsig(\bfu_i)\cdot \nabla  &= \bff_i& \qquad &\text{in $\Omega_i$}\label{eq:elast1}
\\ 
 \jump{\bfsig(\bfu)\cdot \bfn} & = {\bf 0} & \qquad  &\text{on $\Gamma$}\label{eq:elast2}
\\ 
\bfu_i-\bfu_\Gamma &= {\bf 0}& \qquad  &\text{on $\Gamma\cap\partial \Omega_i$}\label{eq:elast3}
\end{alignat}
Here we used the notation
\begin{equation}
\llbracket \bfsig(\bfv)\cdot \bfn \rrbracket := \bfsig(\bfv_1)\cdot \bfn_1+\bfsig(\bfv_2)\cdot \bfn_2
\end{equation}
where $\bfn_i$ denotes the outward pointing normal to $\Omega_i$, $\bfsig(\bfu)$ is the stress 
tensor, $\bfeps\left(\bfu\right) = \left[\varepsilon_{ij}(\bfu)\right]_{i,j=1}^2$ 
is the strain tensor with components
\begin{equation}
\varepsilon_{ij}(\bfu ) = \frac{1}{2}\left( \frac{\partial
  u_i}{\partial x_j}+\frac{\partial u_j}{\partial x_i}\right) 
\end{equation}
$\bfI = \left[\delta_{ij}\right]_{i,j=1}^2$ with $\delta_{ij} =1$
if $i=j$ and $\delta_{ij}= 0$ if $i\neq j$, and $\lambda$ and $\mu$ are the 
Lam\'{e} parameters in plane stress, so in terms of Young's modulus, $\rm{E}$, and Poisson's ratio, $\nu$, we have
\begin{equation}
\lambda = \frac{\rm{E}\,\nu}{1-\nu^2},\quad\mu=\frac{\rm{E}}{2\, (1+\nu)}.
\end{equation}

Now, multiplying (\ref{eq:elast1}) by test functions $\bfv_i$, $\bfv_i = {\bf 0}$ on $\partial\Omega\setminus\Gamma$, integrating by parts over $\Omega_i$, and using (\ref{eq:elast3}) we find
\begin{align}
\sum_i (\bff_i,\bfv_i)_{\Omega_i}  = {}& \sum_i (-\bfsig(\bfu_i)\cdot \nabla ,\bfv_i)_{\Omega_i} \\
={}&\sum_i (\bfsig(\bfu_i),\bfeps(\bfv_i))_{\Omega_i} -\sum_i (\bfsig(\bfu_i)\cdot\bfn_i,\bfv_i)_{\partial\Omega_i}\\
={}&\sum_i (\bfsig(\bfu_i),\bfeps(\bfv_i))_{\Omega_i} -\sum_i (\bfsig(\bfu_i)\cdot\bfn_i,\bfv_i)_{\Gamma}\\
&-\sum_i (\bfu_i-\bfu_\Gamma,\bfsig(\bfv_i)\cdot\bfn_i)_{\Gamma}+(\gamma_i(\bfu_i-\bfu_\Gamma,\bfv_i-\bfv_\Gamma)_{\Gamma}
\end{align}
where $\gamma_i \in \IR_+$ are arbitrary. Using the interface condition (\ref{eq:elast2}), in the form $(\jump{\bfsig(\bfu)\cdot\bfn},\bfu_\Gamma)_{\Gamma}=0$, we finally obtain
\begin{align}
\sum_i (\bff_i,\bfv_i)_{\Omega_i}  = {}& \sum_i (\bfsig(\bfu_i),\bfeps(\bfv_i))_{\Omega_i} -\sum_i (\bfsig(\bfu_i)\cdot\bfn_i,\bfv_i-\bfv_\Gamma)_{\Gamma}\\
&-\sum_i (\bfu_i-\bfu_\Gamma,\bfsig(\bfv_i)\cdot\bfn_i)_{\Gamma}+(\gamma_i(\bfu_i-\bfu_\Gamma,\bfv_i-\bfv_\Gamma)_{\Gamma}
\end{align}
%Next, we add penalty terms to (\ref{energy1}), yielding an augmented Lagrangian
%\begin{align}
%{\frak L}_\text{aug}(\bfu,\bflam) = {}& \sum_i \frac12(\bfsig(\bfu_i),\bfeps(\bfu_i))_{\Omega_i}-\sum_i(\bflam_i,\bfu_i-\bfu_\Gamma)_\Gamma \\
%&+\frac{1}{2} \sum_i (\gamma(\bfu_i-\bfu_\Gamma),\bfu_i  -\bfu_\Gamma)_{\Gamma}-\sum_i(\bff,\bfu_i)_{\Omega_i}\label{energy2}
%\end{align}
%where $\gamma$ is a penalty parameter.
%

To formulate a finite element method we let $V^h_i$, $i=1,2$, and $V^h_\Gamma$ be conforming finite element spaces, such that $\bfu_i^h\in V^h_i$ and $\bfu_\Gamma^h\in V^h_\Gamma$. We then obtain the hybridized Nitsche method from Burman et al.\cite{BuElHaLaLa19}, extended to linear elasticity: find $(\bfu^h_1,\bfu^h_2,\bfu^h_\Gamma)\in V^h_1\oplus V^h_2\oplus V^h_\Gamma$ such that
\begin{align}
\sum_i(\bff_i,\bfv_i)_{\Omega_i} = {}& \sum_i (\bfsig(\bfu^h_i),\bfeps(\bfv_i))_{\Omega_i}-\sum_i(\bfsig(\bfu^h_i)\cdot\bfn_i,\bfv_i-\bfv_\Gamma)_\Gamma \\ 
&-\sum_i(\bfu^h_i-\bfu^h_\Gamma,\bfsig(\bfv_i)\cdot\bfn_i)_\Gamma\\
&+ \sum_i (\gamma_i(\bfu^h_i-\bfu^h_\Gamma),\bfv_i  -\bfv_\Gamma)_{\Gamma}\label{eq:nitschorig}
\end{align}
for all $(\bfv_1,\bfv_2,\bfv_\Gamma)\in V^h_1\oplus V^h_2\oplus V^h_\Gamma$. % with $\gamma_{0,i}$ sufficiently large, we obtain a symmetric, positive definite system of equations. 

\begin{remark}\label{remark1}
With local meshsize $h_i$ on $\Omega_i$ and the choice 
\begin{equation}
\gamma_i = \gamma_{0,i} h_i^{-1}
\end{equation}
it is possible to 
show that the bilinear form is coercive on the finite element space, provided the parameters $\gamma_{0,i}$ are taken large enough,  which together with Galerkin orthogonality and approximation properties of the finite element spaces leads to optimal order a priori error estimates. We refer to Burman et al.\cite{BuElHaLaLa19} for details.
\end{remark}

\section{Interfaces with Bending and Membrane Stiffness}

\subsection{Strong Coupling}

We now add bending and membrane stiffness of the interface to our functional in the vein of model coupling in creeping flow \cite{MaJaRo05}. 
We assume that we are given and arclength parameter $s$ and a unit tangent vector $\bft$ along $\Gamma$, which creates an orthonormal system with the unit normal $\bfn$ to $\Gamma$. For definiteness we assume that $\bfn = \bfn_1 = -\bfn_2$, and $\bft = \bft_1 = -\bft_2$.
We further split $\bfu_\Gamma$ into a normal  and a tangential part 
\begin{equation}
\bfu_\Gamma = u_n \bfn + u_t \bft
\end{equation}
with $u_t=\bft\cdot\bfu_\Gamma$ and $u_n = \bfn \cdot \bfu_\Gamma$.
%Our modified
%augmented functional can then be written
%\begin{align}
%{\frak L}_\text{aug}(\bfu,\bflam) = {}& \sum_i \frac12(\bfsig(\bfu_i),\bfeps(\bfu_i))_{\Omega_i}-\sum_i(\bflam_i,\bfu_i-\bfu_\Gamma)_\Gamma \\ 
%&+\frac{1}{2} \sum_i (\gamma(\bfu_i-\bfu_\Gamma),\bfu_i  -\bfu_\Gamma)_{\Gamma}+\frac12\left(\text{EA}\frac{du_t}{ds},\frac{du_t}{ds}\right)_\Gamma\\ 
%&+\frac12\left(\text{EI}\frac{d^2u_n}{ds^2},\frac{d^2u_n}{ds^2}\right)_\Gamma  -\sum_i(\bff,\bfu_i)_{\Omega_i}
%\end{align}

The equilibrium equations on $\Gamma$ are then assumed as follows
\begin{equation}
\frac{d^2}{ds^2}\left(\text{EI}\frac{d^2u_n}{ds^2}\right)  = f_n-\bfn\cdot\jump{\bfsig\cdot\bfn}
\end{equation}
and
\begin{equation}
-\frac{d}{ds}\left(\text{EA}\frac{du_t}{ds}\right)  = f_t-\bft\cdot\jump{\bfsig\cdot\bfn}
\end{equation}
where $f_n$ and $f_t$ are given external loads. Here \text{EI} denotes bending stiffness (with \text{I} the second moment of inertia) and \text{EA} axial stiffness (with \text{A} the cross section area), both possibly varying with position.
These two equilibrium equations now replace the interface equilibrium (\ref{eq:elast2}) which no longer holds.

Again, multiplying (\ref{eq:elast1}) by test functions $\bfv_i$, $\bfv_i = {\bf 0}$ on $\partial\Omega\setminus\Gamma$, integrating by parts over $\Omega_i$, and using (\ref{eq:elast3}) we find
\begin{align}
\sum_i (\bff_i,\bfv_i)_{\Omega_i}  = {}& \sum_i - ( \bfsig(\bfu_i) \cdot \nabla ,\bfv_i)_{\Omega_i} \\
={}&\sum_i (\bfsig(\bfu_i),\bfeps(\bfv_i))_{\Omega_i} -\sum_i (\bfsig(\bfu_i)\cdot\bfn_i,\bfv_i)_{\partial\Omega_i}\\
={}&\sum_i (\bfsig(\bfu_i),\bfeps(\bfv_i))_{\Omega_i} -\sum_i (\bfsig(\bfu_i)\cdot\bfn_i,\bfv_i)_{\Gamma}\\
&-\sum_i (\bfu_i-\bfu_\Gamma,\bfsig(\bfv_i)\cdot\bfn_i)_{\Gamma}+(\gamma_i (\bfu_i-\bfu_\Gamma),\bfv_i)_{\Gamma}
\end{align}
Writing $\bfv_\Gamma = v_n\bfn + v_t\bft$ we see that
\begin{align}
\sum_i(\bfsig(\bfu_i)\cdot\bfn_i,\bfv_\Gamma)_{\Gamma} 
&= (\jump{\bfsig\cdot\bfn},\bfv_\Gamma)_{\Gamma}
\\
&=\left( \bft\cdot\jump{\bfsig\cdot\bfn},v_t\right)_\Gamma+\left(\bfn\cdot\jump{\bfsig\cdot\bfn},v_n\right)_\Gamma 
\\
&= \left(f_t + \frac{d}{ds}\left(\text{EA}\frac{du_t}{ds}\right),v_t\right)_\Gamma 
\\
&\qquad +\left(f_n-\frac{d^2}{ds^2}\left(\text{EI}\frac{d^2u_n}{ds^2}\right),v_n\right)_\Gamma
\end{align}
%
%Seeking stationary points and using that $(\bfu_i-\bfu_\Gamma)\vert_\Gamma = {\bf 0}$, we then find that
%\begin{align}
%\sum_i(\bff,\bfu_i)_{\Omega_i} ={}&  \sum_i (\bfsig(\bfu_i),\bfeps(\bfv_i))_{\Omega_i}   -\sum_i(\bflam_i,\bfv_i-\bfv_\Gamma)_\Gamma \\ 
%&+\left(\text{EA}\frac{du_t}{ds},\frac{dv_t}{ds}\right)_\Gamma
%+\left(\text{EI}\frac{d^2u_n}{ds^2},\frac{d^2v_n}{ds^2}\right)_\Gamma  
%\end{align}
%Integration by parts and using the equilibrium equations yields
%\begin{align}
% {\bf 0} ={}&  \sum_i(\bfsig(\bfu_i)\cdot\bfn_i,\bfv_i)_{\Gamma}   -\sum_i(\bflam_i,\bfv_i-\bfv_\Gamma)_\Gamma \\ 
%&-\left( \bft\cdot\jump{\bfsig\cdot\bfn},v_t\right)_\Gamma
%-\left(\bfn\cdot\jump{\bfsig\cdot\bfn},v_n\right)_\Gamma  
%\end{align}
%and thus we find
%\begin{equation}
%\sum_i(\bfsig(\bfu_i)\cdot\bfn_i,\bfv_i-\bfv_\Gamma)_{\Gamma}  =\sum_i(\bflam_i,\bfv_i-\bfv_\Gamma)_\Gamma 
%\end{equation}
%and again $\bflam_i= \bfsig(\bfu_i) \cdot \bfn_i$.
Our Nitsche method thus takes the form: find $(\bfu^h_1,\bfu^h_2,u_n^h,u_t^h)\in V^h_1\oplus V^h_2\oplus V^h_n\oplus V^h_t$ such that
\begin{align}
&\sum_i (\bfsig(\bfu^h_i),\bfeps(\bfv_i))_{\Omega_i}  -\sum_i(\bfsig(\bfu^h_i)\cdot\bfn_i,\bfv_i-v_n\bfn-v_t\bft)_\Gamma 
\\ 
&\qquad\qquad   -\sum_i(\bfu^h_i-u_n^h\bfn-u_t^h\bft,\bfsig(\bfv_i)\cdot\bfn_i)_\Gamma
\\
&\qquad \qquad+ \sum_i (\gamma_i(\bfu^h_i-u_n^h\bfn-u_t^h\bft),\bfv_i  -v_n\bfn-v_t\bft)_{\Gamma} 
\\ 
&\qquad \qquad  +\left(\text{EA}\frac{du^h_t}{ds},\frac{dv_t}{ds}\right)_\Gamma+\left(\text{EI}\frac{d^2u^h_n}{ds^2},\frac{d^2v_n}{ds^2}\right)_\Gamma
\\
&\qquad=\sum_i(\bff,\bfv_i)_{\Omega_i}+(f_t,v_t)_{\Gamma}+(f_n,v_n)_{\Gamma}
\end{align}
for all $(\bfv_1,\bfv_2,v_n,v_t)\in V^h_1\oplus V^h_2\oplus V^h_n\oplus V^h_t$.

\begin{itemize}
\item 
In this setting, $V^h_n$ must be a space of $C^1(\Gamma)-$continuous polynomials, whereas $V^h_t$ can be $C^0(\Gamma)-$continuous. Thus it is reasonable to choose different discretizations for $u_t^h$ and $u^h_n$. 
\item 
 For an interface that have a corner or bifurcates in a points we can not use $(u_n^h,u^h_t)$ as global degrees of freedom since these are not continuous if the segments meet at an angle. Thus we must transform the variables 
 back to Cartesian coordinates. We give details for our chosen discretization below.
\end{itemize}

For brevity let us define
\begin{equation}
A(\bfu,\bfv) := \sum_i (\bfsig(\bfu_i),\bfeps(\bfv_i))_{\Omega_i}+\left(\text{EA}\frac{du_t}{ds},\frac{dv_t}{ds}\right)_\Gamma+\left(\text{EI}\frac{d^2u_n}{ds^2},\frac{d^2v_n}{ds^2}\right)_\Gamma
\end{equation}
\begin{equation}
 L(\bfv) :=\sum_i(\bff,\bfv_i)_{\Omega_i}+(f_t,v_t)_{\Gamma}+(f_n,v_n)_{\Gamma}
\end{equation}
The finite element method takes the form: find $\bfu^h :=(\bfu^h_1,\bfu^h_2,u_n^h,u_t^h)\in V^h:= V^h_1\oplus V^h_2\oplus V^h_n\oplus V^h_t$
such that
\begin{align}
L(\bfv) = {}& A(\bfu^h,\bfv)-\sum_i(\bfsig(\bfu^h_i)\cdot\bfn_i,\bfv_i-\bfv_\Gamma)_\Gamma \\
& -\sum_i(\bfu^h_i-\bfu_\Gamma^h,\bfsig(\bfv_i)\cdot\bfn_i)_\Gamma \\
& + \sum_i (\gamma_i(\bfu^h_i-\bfu^h_\Gamma),\bfv_i  -\bfv_\Gamma)_{\Gamma} ,\quad\forall \bfv\in V^h
\end{align}
where we recall that 
\begin{align}
\bfu^h_\Gamma = u_n^h \bfn + u_t^h \bft, \qquad \bfv_\Gamma = v_n \bfn + v_t \bft
\end{align}

\subsection{Cohesive Coupling}

To model a weaker coupling we proceed in the spirit of Juntunen and Stenberg \cite{JuSt09} and Hansbo and Hansbo \cite{HaHa04}.
The cohesive model, replacing the strong condition $\bfu_i=\bfu_\Gamma$ on $\Gamma$, is given by
\begin{equation}\label{eq:stiffnessform}
\bfsig(\bfu_i) \cdot \bfn_i + \bfS_i (\bfu_i-\bfu_\Gamma) = {\bf 0}\quad\text{on}\;\partial\Omega_i\cap\Gamma
\end{equation}
where $\bfS_i$ are coupling stiffness matrices, assumed to be of the form 
\begin{equation}
\bfS_i = \frac{1}{\alpha_i} \bfn\otimes\bfn + \frac{1}{\beta_i} \bft\otimes \bft
\end{equation}
where $1/\alpha_i$ and $1/\beta_i$ are stiffness parameters normal and tangential to the interface, respectively.
We are interested in the case where $\alpha_i$ and $\beta_i$ can be arbitrarily small, so we write (\ref{eq:stiffnessform}) as
\begin{equation}
\bfC_i\bfsig(\bfu_i) \cdot \bfn_i + \bfu_i-\bfu_\Gamma = {\bf 0}
\end{equation}
where $\bfC_i=\bfS_i^{-1} = \alpha_i\bfn\otimes\bfn + \beta_i \bft\otimes \bft$. We then find
\begin{align}
L(\bfv) ={}&A(\bfu,\bfv) - \sum_i(\bfsig(\bfu_i)\cdot\bfn_i,\bfv_i-\bfv_\Gamma)_\Gamma
\\
={}&A(\bfu,\bfv) + \sum_i(\bfsig(\bfu_i)\cdot\bfn_i,\bfC_i\bfsig(\bfv_i)\cdot\bfn_i )_{\Gamma} \\
& -\sum_i(\bfsig(\bfu_i)\cdot\bfn_i,\bfC_i\bfsig(\bfv_i)\cdot\bfn_i+\bfv_i-\bfv_\Gamma)_{\Gamma}\\
={}&A(\bfu,\bfv) + \sum_i(\bfsig(\bfu_i)\cdot\bfn_i,\bfC_i\bfsig(\bfv_i)\cdot\bfn_i )_{\Gamma} \\
& -\sum_i(\bfsig(\bfu_i)\cdot\bfn_i,\bfC_i\bfsig(\bfv_i)\cdot\bfn_i+\bfv_i-\bfv_\Gamma)_{\Gamma}\\
& -\sum_i(\bfC_i\bfsig(\bfu_i)\cdot\bfn_i+\bfu_i-\bfu_\Gamma,\bfsig(\bfv_i)\cdot\bfn_i)_{\Gamma}\\
& + \sum_i(\bfC_i\bfsig(\bfu_i)\cdot\bfn_i+\bfu_i-\bfu_\Gamma , \bftau_i (\bfC_i\bfsig(\bfu_i)\cdot\bfn_i+\bfu_i-\bfu_\Gamma) )_{\Gamma}
\label{eq:cohesive}\end{align}
where the last two terms are zero due to the interface condition and the resulting form on the right hand side is symmetric. Furthermore, $\bftau_i$ is a stabilization matrix of the form
\begin{align}\label{eq:tau-def}
\bftau_i =  \tau_n^i \bfn \otimes \bfn+\tau_t^i \bft \otimes \bft , \qquad \tau_n^i = \frac{1}{  h_i/\gamma_{0 ,i}+ \alpha_i}, \quad 
\tau_t^i = \frac{1}{ h_i/\gamma_{0,i} + \beta_i}
\end{align} 
where $\gamma_{0,i}$ is sufficiently large (cf. Remark \ref{remark1}). 

\subsection{One-Sided Cohesive Coupling}

It is clear that the cohesive model is unphysical in that the beams are allowed to penetrate the domains. Thus we need to 
enforce strong continuity of contact type in such situations.
We then enforce the contact constraints by way of $\bfn_i\cdot(\bfu_i-u_n\bfn)\leq 0$.
With $\sigma_n:= \bfn\cdot\bfsig\cdot\bfn$, $\sigma_t:=\bft\cdot\bfsig\cdot\bfn$, $\jump{v_n^i} := \bfn_i\cdot(\bfv_i-v_n\bfn)$, and
$\jump{v_t^i} := \bft_i\cdot(\bfv_i-v_t\bft)$ our contact conditions on $\Gamma$ can then be formulated, following Burman and Hansbo\cite{BuHa17}: 
\begin{align}
\beta_i\sigma_t(\bfu_i) +\jump{u_t^i}= {}& 0  \quad \text{on}\quad\partial\Omega_i\cap\Gamma\\
\jump{u_n^i}\leq {}& 0\quad \text{on}\quad\partial\Omega_i\cap\Gamma, \label{eq:kuhn1}\\
\sigma_n(\bfu_i) +\alpha_i^{-1}\jump{u_n^i} \leq {}& 0  \quad \text{on}\quad\partial\Omega_i\cap\Gamma\label{eq:kuhn2} \\
(\sigma_n(\bfu_i) +\alpha_i^{-1}\jump{u_n^i} )\jump{u_n^i} = {}& 0  \quad \text{on}\quad\partial\Omega_i\cap\Gamma\label{eq:kuhn3}
\end{align}
where we recognise (\ref{eq:kuhn1})--(\ref{eq:kuhn3}) as the Kuhn--Tucker conditions. We begin by rewriting (\ref{eq:cohesive})
in normal and tangential components, writing :
\begin{align}
L(\bfv) 
&=A(\bfu,\bfv) + \sum_i(\beta_i\sigma_t(\bfu_i),\sigma_t(\bfv_i))_{\Gamma} 
\\
&\qquad
-\sum_i(\sigma_t(\bfu_i),\beta_i\sigma_t(\bfv_i)+\jump{v_t^i})_{\Gamma}
\\
&\qquad  -\sum_i(\beta_i\sigma_t(\bfu_i)+\jump{u_t^i},\sigma_t(\bfv_i))_{\Gamma} 
\\
&\qquad  + \sum_i(\beta_i\sigma_t(\bfv_i)+\jump{v_t^i} ,\tau_t^i( \beta_i\sigma_t(\bfu_i)+\jump{u_t^i}) )_{\Gamma}
\\
&\qquad -\sum_i(\sigma_n(\bfu_i),\jump{v_n^i})_{\Gamma}
\end{align}
where we did not introduce the normal component of the cohesive law.

We now turn to the alternative formulation of the Kuhn--Tucker conditions due to Rockafellar \cite{Rock74}, introduced in a Nitsche formulation for contact 
analysis by Chouly and Hild \cite{ChHi13}. We recognize that 
\begin{equation}
p_i :=\sigma_n(\bfu_i) +\alpha_i^{-1}\jump{u_n^i}
\end{equation}
act as multipliers (cf. Burman and Hansbo\cite{BuHa17}, Section 5.1), and
the Kuhn--Tucker conditions are then equivalent to the relation
\begin{equation}
p_i= -\frac{1}{\epsilon_i}[ \jump{u_n^i} -\epsilon_ip_i]_+
\end{equation}
where $[x]_+ = \max(x,0)$ and $\epsilon_i > 0$ but arbitrary.
We then write
\begin{align}
(\sigma_n(\bfu_i),\jump{v_n^i})_{\Gamma}=  {}&(\sigma_n(\bfu_i)+\alpha_i^{-1}\jump{u_n^i},\jump{v_n^i})_{\Gamma} -(\alpha_i^{-1}\jump{u_n^i},\jump{v_n^i})_{\Gamma}
  \\
 =  {}&(\sigma_n(\bfu_i)+\alpha_i^{-1}\jump{u_n^i},\jump{v_n^i}-\epsilon_i(\sigma_n(\bfv_i)+\alpha_i^{-1}\jump{v_n^i}))_{\Gamma}\\
& +(\epsilon_i(\sigma_n(\bfu_i)+\alpha_i^{-1}\jump{u_n^i}),\sigma_n(\bfv_i)+\alpha_i^{-1}\jump{v_n^i})_{\Gamma}\\
& -(\alpha_i^{-1}\jump{u_n^i},\jump{v_n^i})_{\Gamma}
   \end{align}
Thus we have that
\begin{align}
& -(\sigma_n(\bfu_i),\jump{v_n^i})_{\Gamma}=(\alpha_i^{-1}\jump{u_n^i},\jump{v_n^i})_{\Gamma} \nonumber \\
&+\left(\epsilon_i^{-1}\left[\jump{u_n^i}-\epsilon_i(\sigma_n(\bfu_i)+\alpha_i^{-1}\jump{u_n^i})\right]_+,\jump{v_n^i}-\epsilon_i(\sigma_n(\bfv_i)+\alpha_i^{-1}\jump{v_n^i})\right)_{\Gamma} \nonumber \\ \nonumber 
& -(\epsilon_i(\sigma_n(\bfu_i)+\alpha_i^{-1}\jump{u_n^i}),\sigma_n(\bfv_i)+\alpha_i^{-1}\jump{v_n^i})_{\Gamma} 
\end{align}
To obtain well conditioned systems for small $\alpha_i$, we examine the case when $p_i = 0$ and $p_i\neq 0$.

\begin{itemize}
\item If $p_i\neq 0$ we are in contact and
\begin{align*}
 -(\sigma_n(\bfu_i),\jump{v_n^i})_{\Gamma}={}&-(\alpha_i^{-1}\jump{u_n^i},\jump{v_n^i})_{\Gamma}+(\epsilon_i^{-1}\jump{u_n^i},\jump{v_n^i})_{\Gamma}\\
&-\left(\sigma_n(\bfu_i),\jump{v_n^i}\right)_{\Gamma}
-(\jump{u_n^i},\sigma_n(\bfv_i))_{\Gamma} 
\end{align*}
which, with $\epsilon_i^{-1}=\gamma_{0,i}/h_i+\alpha_i^{-1}$ gives a scheme of the type (\ref{eq:nitschorig}).
\item If $p_i= 0$ we have the cohesive law and
\begin{align*}
 -(\sigma_n(\bfu_i),\jump{v_n^i})_{\Gamma}={}&(\alpha_i^{-1}\jump{u_n^i},\jump{v_n^i})_{\Gamma}\\
&-(\epsilon_i(\sigma_n(\bfu_i)+\alpha_i^{-1}\jump{u_n^i},\sigma_n(\bfv_i)+\alpha_i^{-1}\jump{v_n^i})_{\Gamma}\\
={}&-(\alpha_i\sigma_n(\bfu_i)+\jump{u_n^i},\sigma_n(\bfv_i))_{\Gamma}\\
&-(\sigma_n(\bfu_i),\alpha_i\sigma_n(\bfv_i)+\jump{v_n^i})_{\Gamma}+(\alpha_i\sigma_n(\bfu_i),\sigma_n(\bfv_i))_{\Gamma}\\
&+((\alpha_i-\epsilon_i)(\sigma_n(\bfu_i)+\alpha_i^{-1}\jump{u_n^i}),\sigma_n(\bfv_i)+\alpha_i^{-1}\jump{v_n^i})_{\Gamma}\end{align*}
and using again $\epsilon_i^{-1}=\gamma_{0,i}/h_i+\alpha_i^{-1}$ we find the standard method for cohesive laws with penalty parameter
$(\alpha_i-\epsilon_i)/\alpha_i^2 = 1/(h_i/\gamma_{0,i}+\alpha_i)=\tau_n^i$.
\end{itemize}
Thus the parameter $\epsilon_i$ does not change in contact and out of contact and the schemes become insensitive to small $\alpha_i$.

\section{Numerical examples} 

In this Section we illustrate the properties of the model and method by presenting 
some basic numerical examples. In all cases we used Young's modulus $\text{E} =10^6$ and Poisson's ratio $\nu=1/3$ in the bulk, and $\gamma_{0,i}=20(\lambda+\mu)$ as a Nitsche penalty parameter.

\begin{figure}[ht]
	\begin{center}
		\includegraphics[scale=0.35]{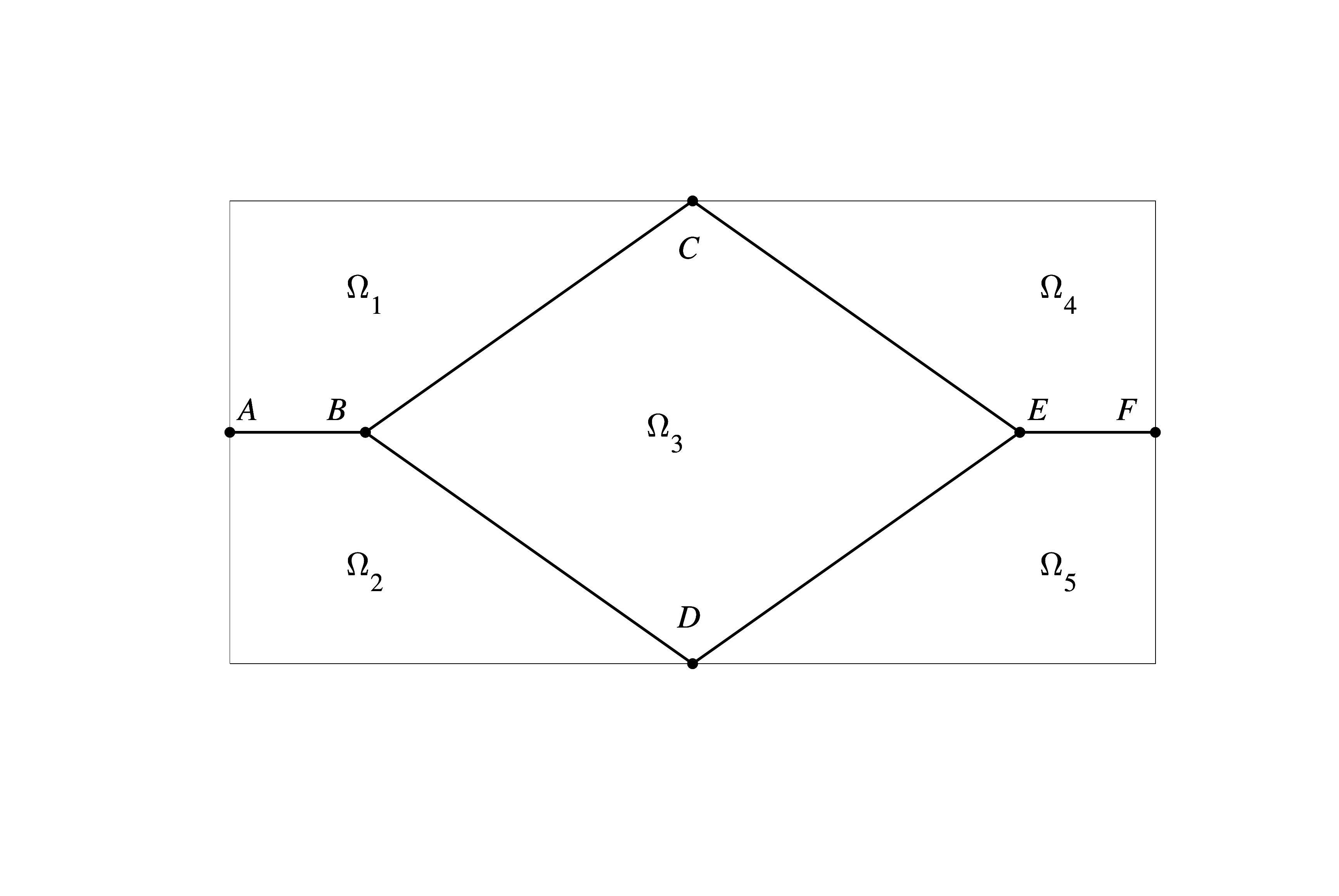}
	\end{center}
	\caption{Cantilever domain.}
	\label{fig:domain}
\end{figure}

We consider a macro domain $\Omega=(0,2)\times (0,1)$ split into 5 subdomains $\Omega_1$ to $\Omega_5$ separated by 6 line segments
$\vv{AB}$, $\vv{BC}$, $\vv{BD}$, $\vv{CE}$, $\vv{DE}$, and $\vv{EF}$, where $A=(0,1/2)$, $B=(1-1/\sqrt{2},1/2)$, $C=(1,1)$, $D=(1,0)$,
$E=(1+1/\sqrt{2},1/2)$, $F=(2,1/2)$, see Fig. \ref{fig:domain}.

We use one cubic $C^1$ polynomial in each segment for the approximation of the interface normal displacement $u_n$, and one linear polynomial in each segment for the tangential displacement $u_t$. Continuity at the endpoints of the segments is achieved by a transformation to Cartesian coordinates of the nodal variables by use of 
\begin{equation}
\left[\begin{array}{>{\displaystyle}c}u_x\\u_y\\ \theta\end{array}\right] =\left[\begin{array}{>{\displaystyle}c>{\displaystyle}c>{\displaystyle}c}n_x & t_x & 0\\ n_y & t_y & 0 \\ 0 & 0 & 1\end{array}\right]\left[\begin{array}{>{\displaystyle}c}u_n\\u_t\\ \theta\end{array}\right]
\end{equation}
where $\theta := u_n'(s)$ is the ``rotation'' degree of freedom in the $C^1$ approximation.

The boundary conditions are: Dirichlet boundary conditions $(u_x,u_y) =(0,0)$ at $x=0$ and zero Neumann conditions elsewhere.
These are also imposed on the interface variables. We further impose zero rotation for the interface variable $u_n$ at $x=0$ in the case $\text{EI} > 0$. The approximation in the domains is a $P^1$--$C^0$ approximation (constant strain triangle). In the case of cohesion, we use
the same constants $\alpha_i =: \alpha$ and $\beta_i =:\beta$ for all interfaces.

\subsection{Bending of a Cantilever Structure} 

\begin{figure}[ht]
	\begin{center}
		\includegraphics[scale=0.25]{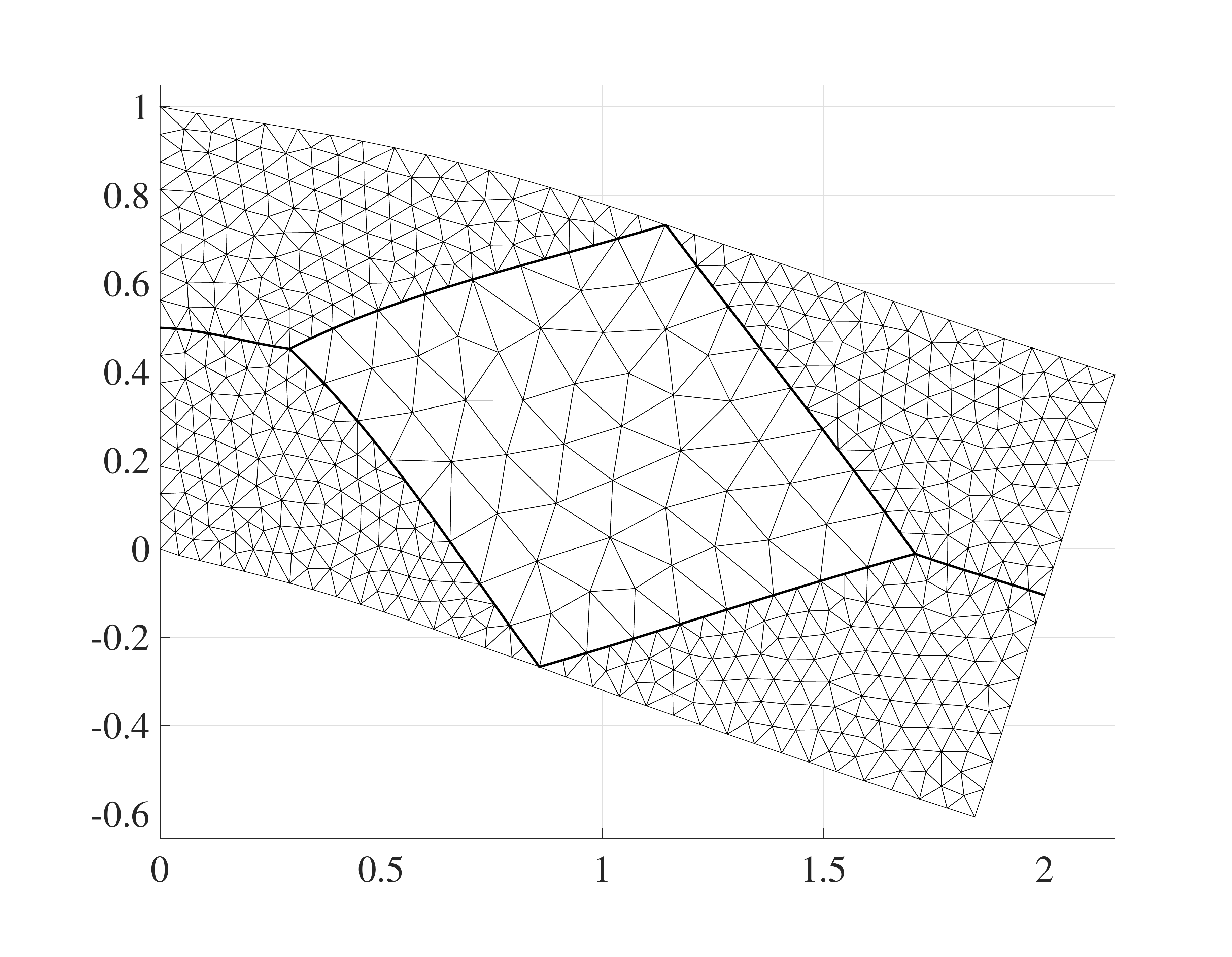}
	\end{center}
	\caption{Deformations for the hybrid method without interface stiffness.}
	\label{fig:hybrid}
\end{figure}

\begin{figure}[ht]
	\begin{center}
		\includegraphics[scale=0.2]{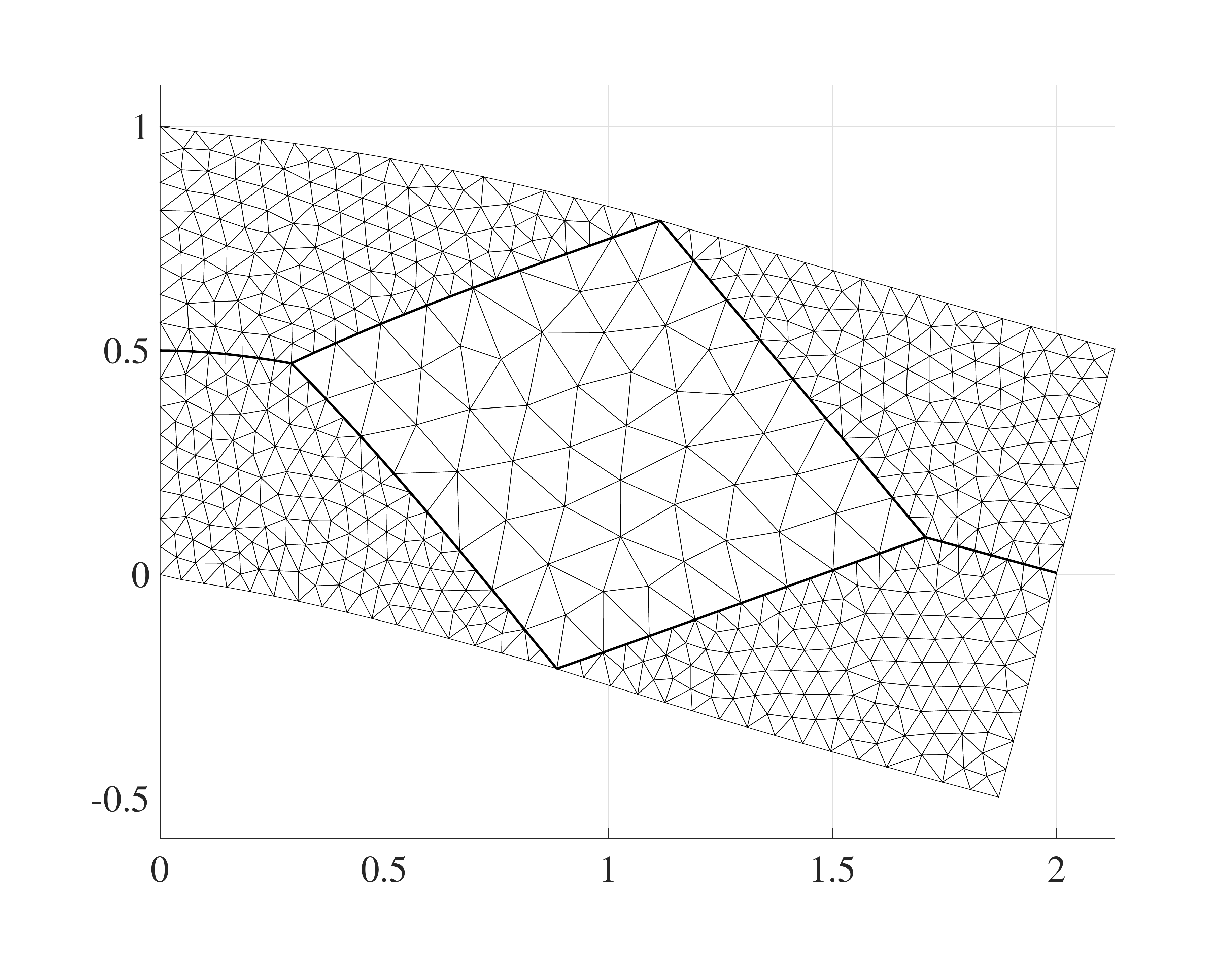}\includegraphics[scale=0.2]{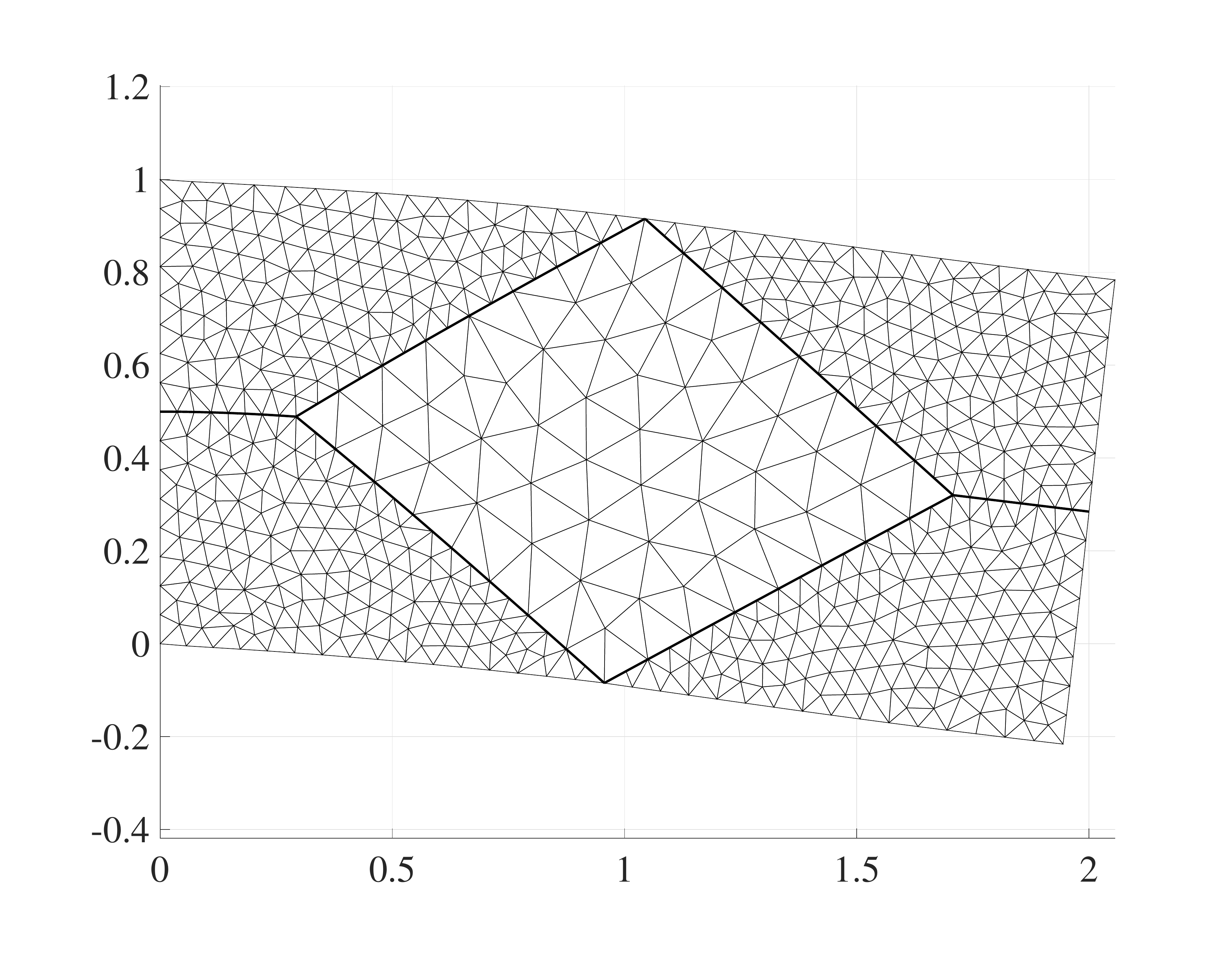}
	\end{center}
	\caption{Deformations with interface bending stiffness, $\text{EI}=10^4$ (left) and $\text{EI}=10^5$ (right).}
	\label{fig:bend}
\end{figure}

We consider constant loads $\bff_i=(0,-2\times 10^4)$, $i=1,\ldots,5$, and show the computational results for the standard hybrid method ($\text{EI}=\text{EA}=\alpha=\beta=0$) in Fig. \ref{fig:hybrid}. We next show the effect of increasing bending stiffness on the interface, with $\text{EI} = 10^4$
and $\text{EI}=10^5$ in Fig. \ref{fig:bend}. The stiffening effect is noticeable.

\begin{figure}[ht]
	\begin{center}
		\includegraphics[scale=0.2]{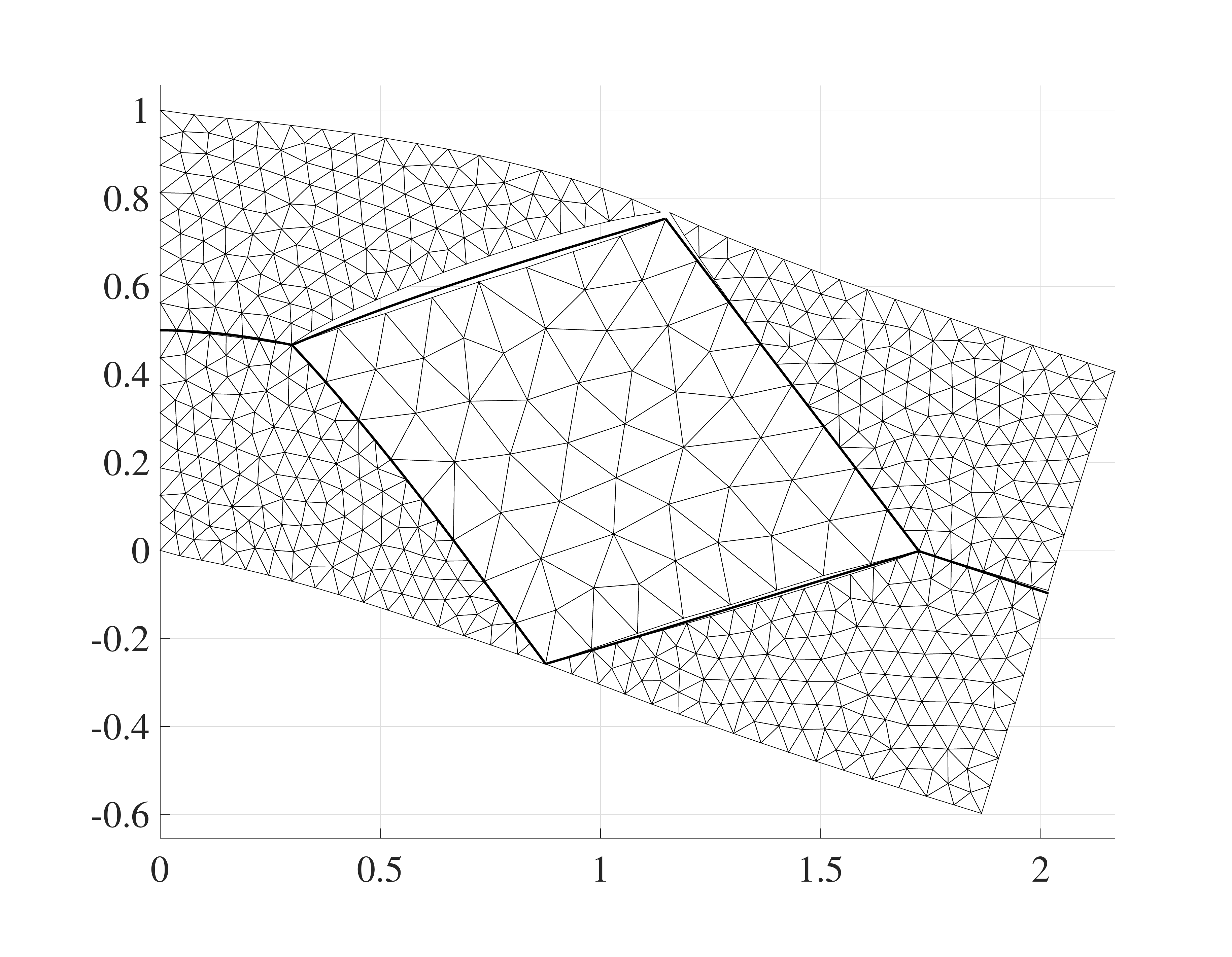}\includegraphics[scale=0.2]{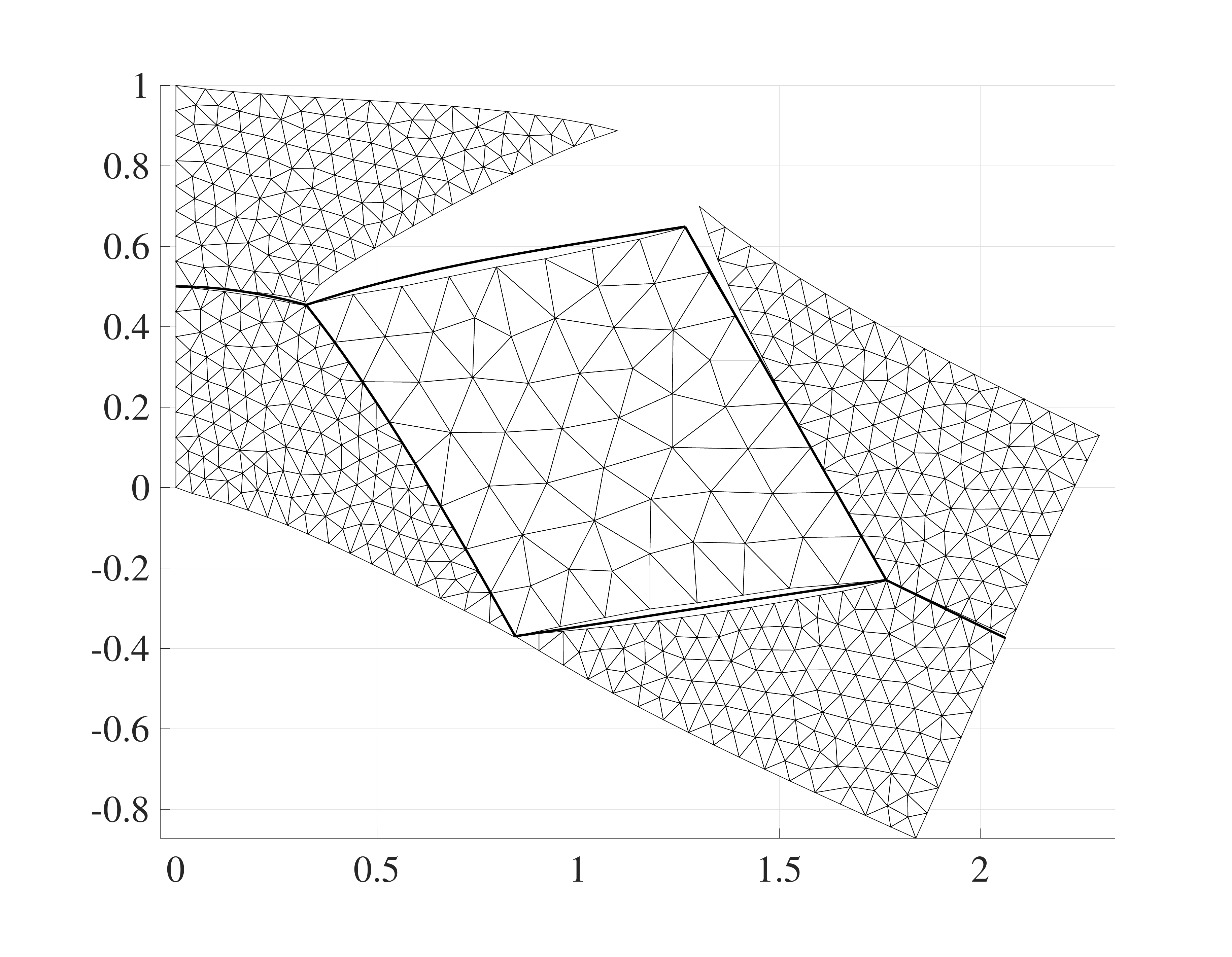}
	\end{center}
	\caption{Deformations with interface bending stiffness and normal cohesion, $\text{EI}=10^4$; $\alpha=10^{-6}$ (left) and $\alpha=10^{-5}$ (right).}
	\label{fig:normal}
\end{figure}

Finally, in Fig. \ref{fig:normal} we show the effect of normal compliance at a fixed bending stiffness $\text{EI}=10^4$. We show the
results for $\alpha=10^{-6}$ and $\alpha=10^{-5}$. The contact algorithm is invoked to avoid domain penetrations.

\subsection{Stretching} 

\begin{figure}[ht]
	\begin{center}
		\includegraphics[scale=0.3]{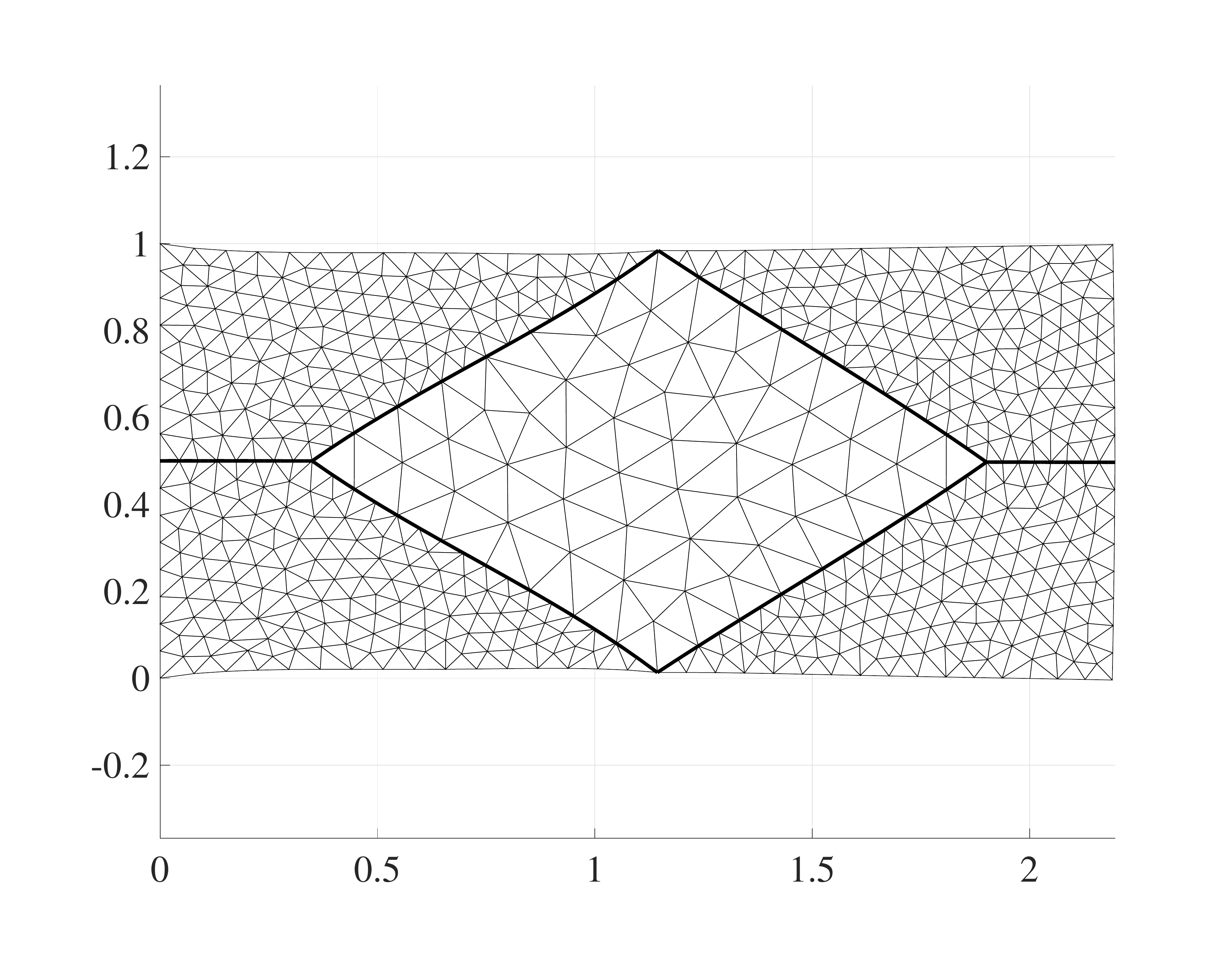}
	\end{center}
	\caption{Stretch deformations for the hybrid method without interface stiffness.}
	\label{fig:hybrids}
\end{figure}

\begin{figure}[ht]
	\begin{center}
		\includegraphics[scale=0.2]{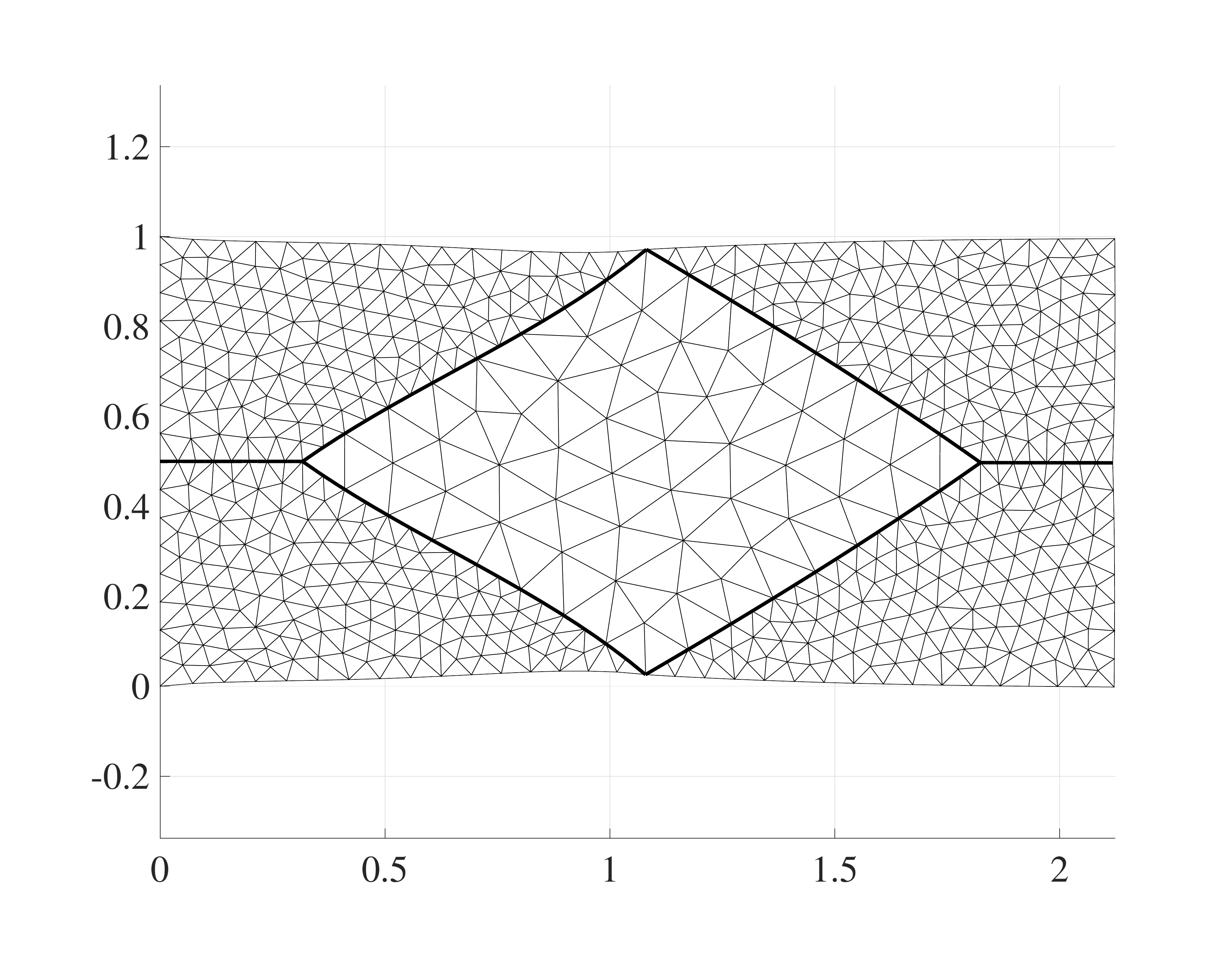}\includegraphics[scale=0.2]{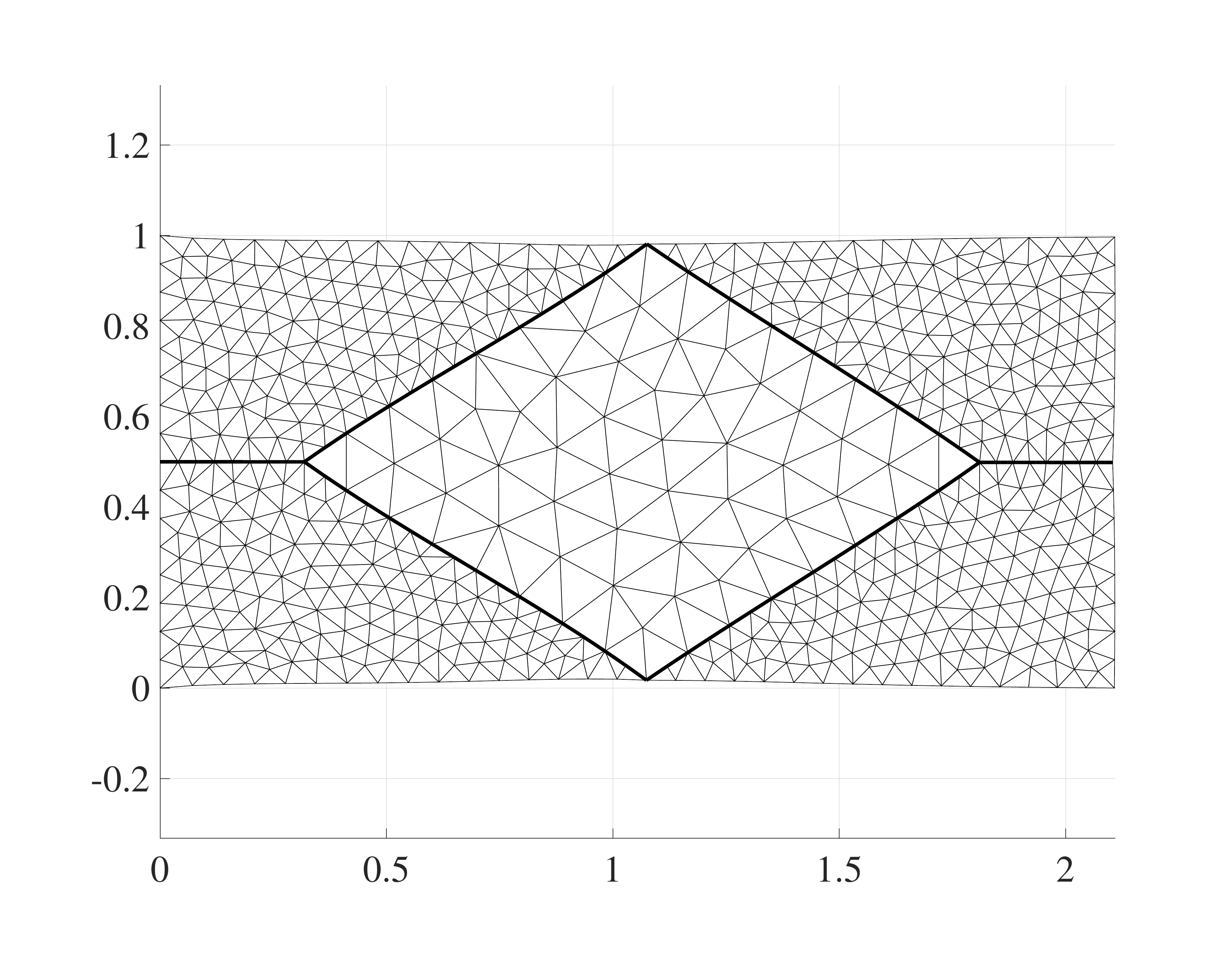}
	\end{center}
	\caption{Deformations with interface membrane stiffness and with combined membrane/bending stiffness, $\text{EA}=10^6$; $\text{EI}=0$ (left) and $\text{EI}=10^{4}$ (right).}
	\label{fig:memb}
\end{figure}

\begin{figure}[ht]
	\begin{center}
		\includegraphics[scale=0.2]{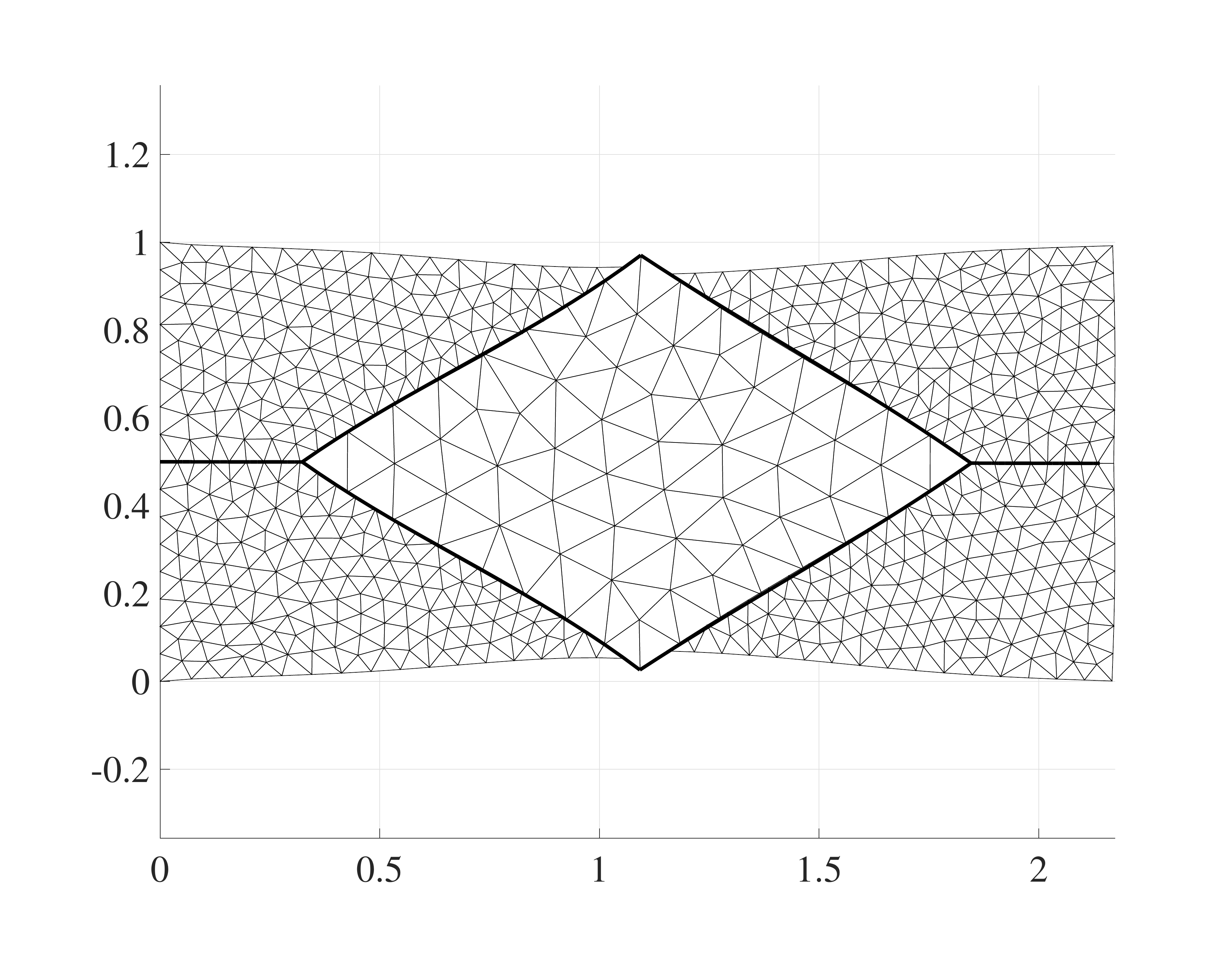}\includegraphics[scale=0.2]{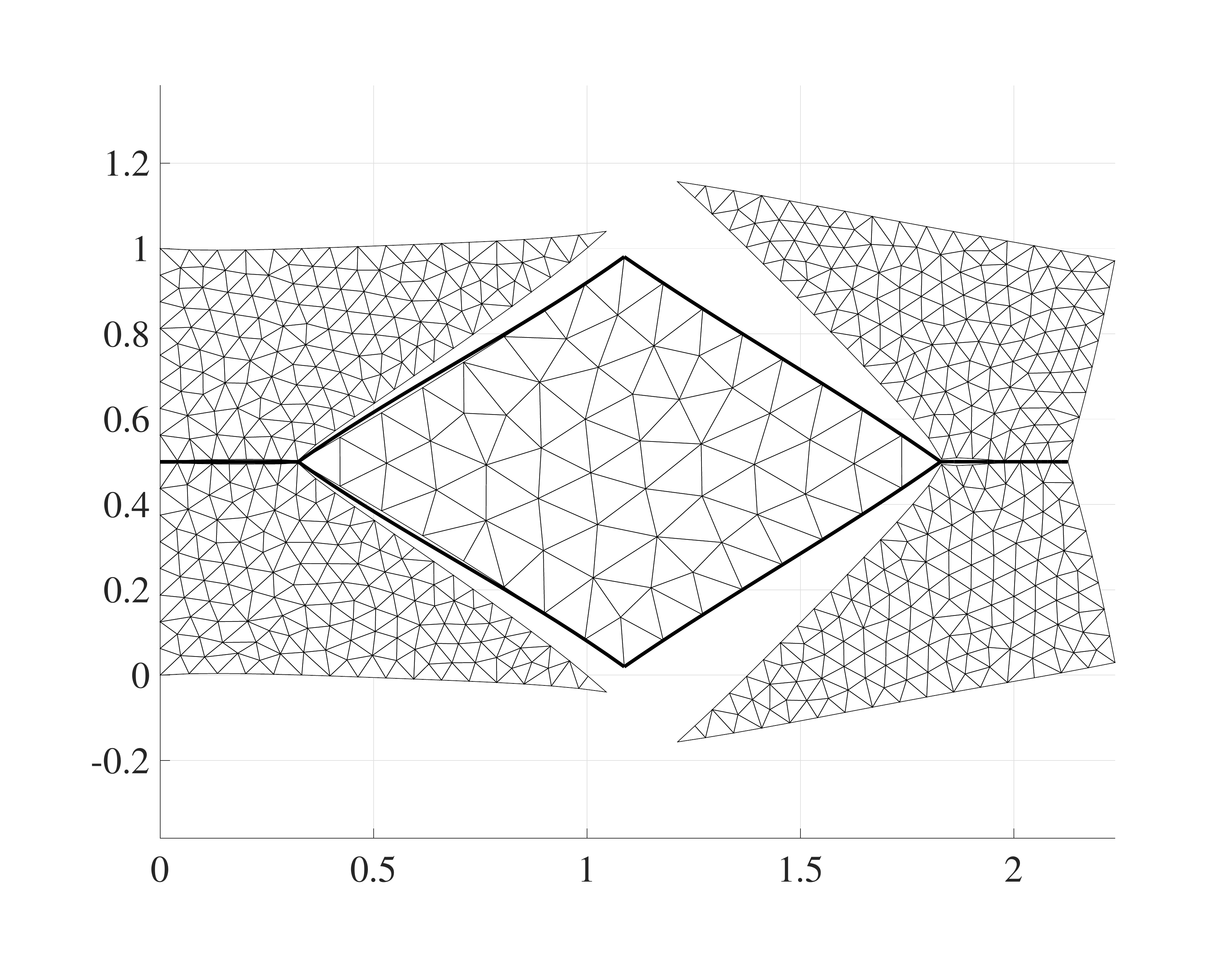}
	\end{center}
	\caption{Deformations with interface stiffness and with cohesion, $\text{EA}=10^6$, $\text{EI}=10^4$; $\beta=10^{-5}$, $\alpha=0$ (left) and $\beta=0$, $\alpha=10^{-5}$ (right).}
	\label{fig:tannorm}
\end{figure}

\begin{figure}[ht]
	\begin{center}
		\includegraphics[scale=0.3]{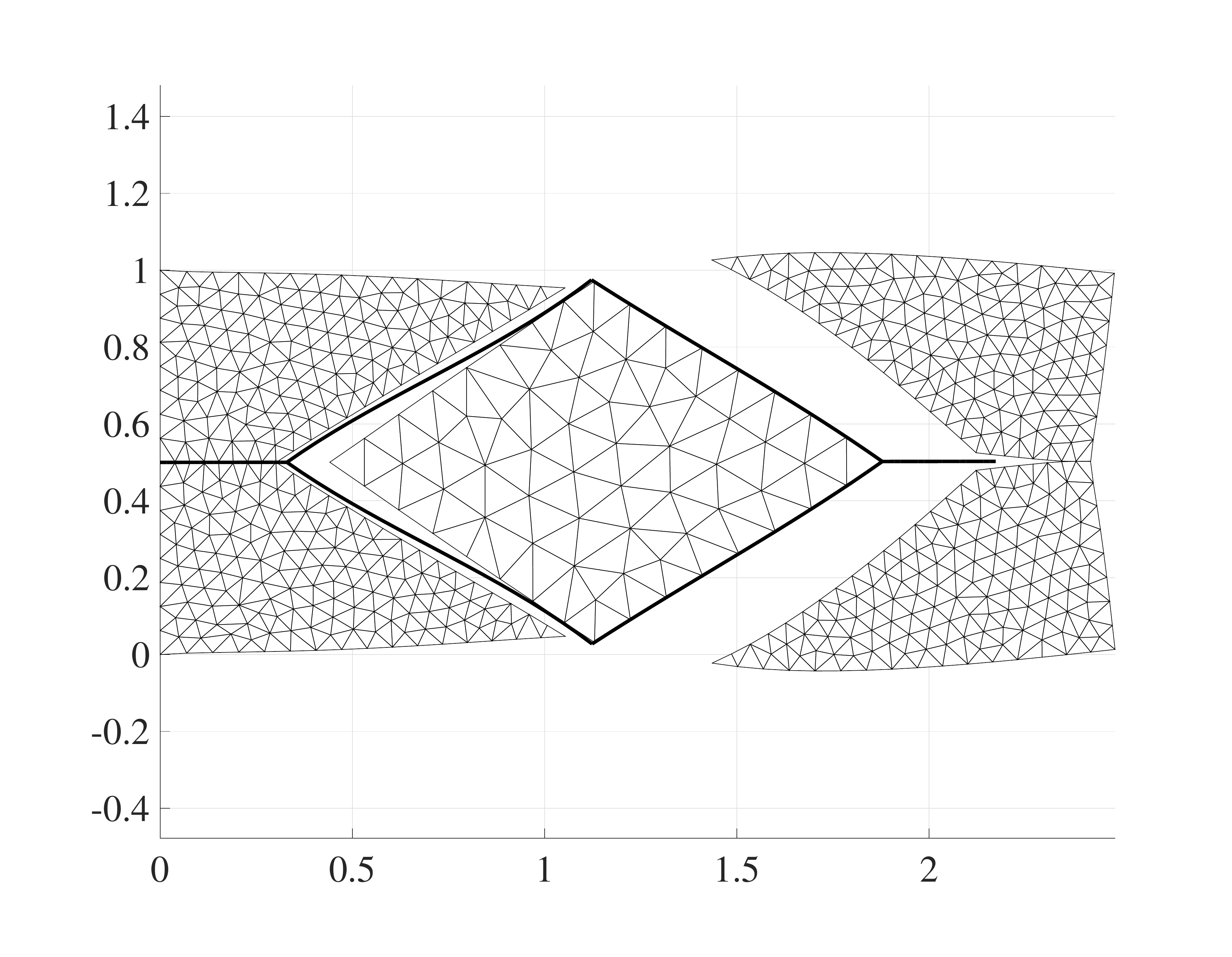}
	\end{center}
	\caption{Deformations with interface stiffness and with cohesion, $\text{EA}=10^6$, $\text{EI}=10^4$; $\beta=10^{-5}$, $\alpha=10^{-5}$.}
	\label{fig:toto}
\end{figure}

The loads in this example are $\bff_i=(10^5,0)$, $i=1,\ldots,5$, to give a stretch of the domain. In Fig. \ref{fig:hybrids} we show the result for the standard hybridized method. We next show the effect of adding membrane and bending stiffness on the interface, with $\text{EA} = 10^6$
and $\text{EI}=0$ and with $\text{EA}=10^6$, $\text{EI}=10^4$ in Fig. \ref{fig:memb}. We then compare the effect of tangential cohesion, $\beta=10^{-5}$, with that of normal cohesion, $\alpha=10^{-5}$, in Fig. \ref{fig:tannorm}. Finally, in Fig. \ref{fig:toto} we show the effect of having both tangential and normal cohesion.

\section{Concluding Remarks} 

We have introduced a hybridized Nitsche method for linearized elasticity which uses an auxiliary interface displacement, modelled independently of the domains. This allows for easy modeling of different stiffness models at the interface. We have focused 
here on Euler--Bernoulli beam bending and membrane stiffness, but other models can be easily accommodated; the only requirement is
that continuity of displacements between the interface field and the domain fields can be represented. We have also 
suggested weaker couplings between the interface and domains in the form of a cohesive interface law, with no-penetration fulfilled.
This leads to a nonlinear contact problem which fits straighforwardly in the general framework of Nitsche's method. Some numerical examples 
are provided to show how different parameter choices affect the solution in bending and in stretching of a plane elasticity problem.

\subsection*{Acknowledgements}
This research was supported in part by
the the Swedish Research Council Grants No. 2017-03911, 2018-05262, and 
the Swedish strategic research programme eSSENCE.

\bibstyle{WileyNJD-AMA}

\newpage

%\vfill
%\bigskip
%\bigskip
%\noindent
%\footnotesize {\bf Authors' addresses:}
%
%\smallskip
%\noindent
%Erik Burman,  \quad \hfill \addressuclshort\\
%{\tt e.burman@ucl.ac.uk}
%
%\smallskip-eps-converted-to.pdf
%\noindent
%Peter Hansbo,  \quad \hfill \addressjushort\\
%{\tt peter.hansbo@ju.se}
%
%\smallskip
%\noindent
%Mats G. Larson,  \quad \hfill \addressumushort\\
%{\tt mats.larson@umu.se}

\end{document}